\documentstyle[11pt]{amsart}
\newcommand{\nz}{\mbox{${\Bbb Z}$}}

\newcommand{\vep}{\varepsilon}
\newcommand{\beq}{\begin{equation}}
\newcommand{\eeq}{\end{equation}}
\newcommand{\bn}{\begin{align}}
\newcommand{\ed}{\end{align}}
\newcommand{\vphi}{\varphi}
\newcommand{\UqgD}{\mbox{$U_q^D(\hat{{\frak g}})\,$}}
\newcommand{\UqdvaD}{\mbox{$U_q^D(\hat{{\frak sl}}_2)\,$}}
\newcommand{\UqnD}{\mbox{$U_q^D(\hat{{\frak sl}}_n)\,$}}
\newcommand{\Uqdva}{\mbox{$U_q(\hat{{\frak sl}}_2)\,$}}
\newcommand{\Uqn}{\mbox{$U_q(\hat{{\frak sl}}_n)\,$}}
\newcommand{\Uqg}{\mbox{$U_q(\hat{{\frak g}})\,$}}
\theoremstyle{plain}
 \newtheorem{thm}{Theorem}
 \newtheorem{prop}[thm]{Proposition}

\theoremstyle{definition}

 {}
\theoremstyle{remark}
\newtheorem{rem}{Remark}

\begin{document}
\begin{center}
\hfill ITEP-TH-22/98\\
\hfill RIMS-1199\\
\bigskip
\end{center}
\title{On the FRTS
approach to quantized current algebras }

\author{Jintai  Ding}
\author{Sergei  Khoroshkin }
\address{Jintai Ding, RIMS, Kyoto University,
 Kyoto 606, Japan}
\address{Sergei Khoroshkin, ITEP,
117259 Moscow, Russia}
\begin{abstract}
We study the possibility to establish $L$-operator's formalism by
Faddeev-Reshetikhin-Takhtajan-Semenov-Tian-Shansky (FRST) for quantized
 current algebras, that is, for quantum affine algebras in the
''new realization '' by V. Drinfeld with the corresponding 
Hopf algebra  structure and for their Yangian counterpart.
 We establish this formalism using the twisting procedure by Tolstoy
 and the second author and explain the problems which  FRST approach
 encounter for quantized current algebras.
We show also that, for the case of
$U_q(\hat {\frak sl}_n)$,
   entries of the L-operators of FRTS type give  the
Drinfeld current operators for the non-simple roots, which we discovered
recently. As an application we deduce  the commutation relations between
these current operators for $U_q(\hat {\frak sl}_3)$.
 \end{abstract}
\maketitle
1991 Math Subject Classification(s): 17B37
\pagestyle{plain}
\setcounter{section}{-1}
\section{Introduction.}

The current realization
 of the quantum affine algebra $U_q(\hat{ \frak{g}})$ and of
the Yangian
was obtained by
Drinfeld \cite{Dr1}, which
he called  ''new realization'' of quantum affine
 algebras and of Yangians. The ''new realization'' came with its own Hopf
 algebra structure. In order to distinguish these Hopf algebras from the
 quantum affine algebras in Drinfeld-Jimbo formulation of quantized
 enveloppping Kac--Moody algebras, we use the name quantized current
 algebras and in the notations mark them with upper index $D$.
  These Hopf algebras
 and their generalizations are now extensively studied \cite{DI1},
 \cite{DI2}, \cite{EF1}, \cite{EF2}.

In the paper \cite{DK}, we extended Drinfeld quantization of current algebras
 in several directions: we constructed current operators for nonsimple roots
 of ${\frak{g}}$, defined the new braid group action in terms of currents
 and presented a description of the universal $R$-matrix in two
 equivalent forms: in a form of
 an infinite product and in a form of certain integrals
 over the current operators.

The aim of this paper is , in addition to \cite{DK}, to develop
Faddeev-Reshetikhin-Takhtajan-Semenov-Tian-Shansky (FRST) approach, where
 the generators are gathered into $L$-operators which satisfy famous ''$RLL$''
 relations, for quantized current algebras. Such a description of quantized
 envelopping algebras has  at least two nice features: first, $L$-operators
 are group-like elements, which means that their entries obey simple
comultiplication rules,  and second, the Cartan-Weyl generators (for
 a specific ordering of roots) can be extracted from the $L$-operators.

In FRST approach, the $L$-operators are given by the projections of the
 tensor components of the universal $R$-matrix to certain finite-dimensional
 representations and all their properties follow immediately from the
 properties of the universal $R$-matrix. Unfortunately, this simple
 prescription does not work directly for quantized current algebras.
The reason is as
 follows. Quantized current algebras are by definition topological algebras
 (otherwise the comultiplication structure, written in Laurent series, is not
 well defined), which act on the tensor category of highest weight
 representations. To the contrary, in the definition of tensor category of
 finite- dimensional
 modules,  we come to problems of  divergence.

However, the vector representation of \UqnD is
 well defined and we can project to it the components of the universal
 $R$-matrix. The corresponding $L$-operators look nice: they are triangular
 and Drinfeld current operators can be easily extracted from them.
  The image of the universal $R$-matrix in tensor product of vector
  is also triangular with $\delta$-function term
 out of diagonal. But it does not satisfy the Yang-Baxter equation.

In order to derive proper commutation relations between the entries of the
 $L$-operators, we use the results of \cite{KT1}, where the Hopf structure
 of $U_q^D(\hat{ \frak{g}})$ was obtained from the Hopf structure in
 $U_q(\hat{ \frak{g}})$  by two equivalent ways:
 either via twisting by certain factor of the universal $R$-matrix of
$U_q(\hat{ \frak{g}})$ or as a limit of twists by  automorphisms of
 affine shifts in the ($q$)-Weyl group of $U_q(\hat{ \frak{g}})$.
 We show that the commutation relations for Gauss coordinates of the
 $L$-operators
for finitely twisted $U_q(\hat{ \frak{g}})$ have the limit in the topology
 of $U_q^D(\hat{ \frak{g}})$ and this limit is the defining relations
 for the currents of  $U_q^D(\hat{ \frak{g}})$. This gives a way to restore
 the correct commutation relations for Gauss coordinates of the $L$-operators
by  looking back to the finite twist.
 We also prove that, if we multiply before the twist
the initial $R$ matrix by the function
  $1-u/v$, then the limit of the corresponding  Yang-Baxter equations
 turns to the relations on $L$-operators with diagonal $R$-matrix, which are
 well defined and are the corollaries of the defining relations for the
 currents of $U_q^D(\hat{ \frak{g}})$.

  The limiting procedure  enables us to  prove that
the entries of the
$L^{\pm}(z)$ produce the current operators for the non-simple roots
discovered in  \cite{DK}; so  we can use the commutation
 relations for the entries of $L$-operators in order to deduce
 the commutation relations the current operators for non-simple roots of
${ \frak{g}}$. We present the results of calculations for
$U_q^D(\hat{{\frak sl}}_3)$.

 All the  arguments above
are valid as well for other types of quantized current
 algebras, e.g. for the Yangian type algebras and for elliptic algebras
 of dynamical type. We can also see from the ideology of twists by
 affine shifts,
 that in general the existence of Drinfeld current realization is in a sense
 equivalent to existence of nontrivial family of automorphisms of the
 basic solution of the Yang-Baxter equation. Only Baxter's $8$-vertex
 $R$-matrix does not have evident symmetries of such a type.

The exposition goes as follows.
We first study the finite twisting coming
from the Weyl group elements  and derive the
image  of the twisted R-matrix of  $U_q{\hat {\frak sl}_n}$ on the
fundamental representation $\Bbb C^n$.
We check that they indeed satisfy Yang-Baxter equation.
 Taking certain limit, compatible with a topology
of \UqnD,
we derive  the R-matrix with
entries of singular function $\delta (z)$ and triangular
property.  Then we will study the related quasi-triangular
structure and explain the problem to use
 this R-matrix to formulate
certain  FRTS   construction of quantized current algebras.
 Finally we proceed to derive  the commutation
relation between the Gauss coordinates and show their connection with
the current operators for non-simple roots.
We end with a  discussion of  other types of quintized current algebras.
\section{ Drinfeld realization as twisted quantum affine algebra.}
This section is mainly to  introduce the results from \cite{KT1}.
We will use basically   the notation from \cite{KT1}, which
we refer the reader to.
In \cite{KT1}, an  idea form Drinfeld \cite{Dr3} and Reshetikhin \cite{R},
where Drinfeld studies the twisting of the Hopf algebra structure is
utilized.

Let $A=(a_{ij})$ $i,j=1,...,r$ be the Cartan matrix of simple Lie
algebra ${ \frak{g}}$
 of simple laced type. We define the quantized current algebra
$U_q^D(\hat{ \frak{g}})$ as follows.

The Hopf algebra $U_q^D(\hat{ \frak{g}})$ is an associative algebra with unit
1 and the generators: $\vphi_i(m)$,$\psi_i(-m)$, $x^{\pm}_i(l)$, for
$i=i,...,r$, $l\in \nz
$ and $m\in -\nz_+$ and a central
 element $c$. Let $z$ be a formal variable and
 $x_i^{\pm}(z)=\sum_{l\in \nz}x_i^{\pm}(l)z^{-l}$,
$\vphi_i(z)=\sum_{m\in -\nz_+}\vphi_i(m)z^{-m}$ and
$\psi_i(z)=\sum_{m\in \nz_+}\psi_i(m)z^{-m}$. In terms of the
formal variables,
the defining relations are
\begin{align*}
& \vphi_i(z)\vphi_j(w)=\vphi_j(w)\vphi_i(z), \\
& \psi_i(z)\psi_j(w)=\psi_j(w)\psi_i(z),\\
& \vphi_i(z)\psi_j(w)\vphi_i(z)^{-1}\psi_j(w)^{-1}=
  \frac{g_{ij}(\frac z wq^{-c})}{g_{ij}(\frac z wq^{c})}, \\
& \vphi_i(z)x_j^{\pm}(w)\vphi_i(z)^{-1}=
  g_{ij}(\frac z wq^{\mp \frac{1}{2}c})^{\pm1}x_j^{\pm}(w),
\end{align*}
\begin{align*}
& \psi_i(z)x_j^{\pm}(w)\psi_i(z)^{-1}=
  g_{ij}(\frac w zq^{\mp \frac{1}{2}c})^{\mp1}x_j^{\pm}(w), \\
& [x_i^+(z),x_j^-(w)]=\frac{\delta_{i,j}}{q-q^{-1}}
  \left\{ \delta(\frac z wq^{-c})\psi_i(wq^{\frac{1}{2}c})-
          \delta(\frac z wq^{c})\vphi_i(zq^{\frac{1}{2}c}) \right\}, \\
& (zq^{}-q^{\pm a_{ij}}w)x_i^{\pm}(z)x_j^{\pm}(w)=
  (q^{\pm a_{ij}}z-wq^{-j})x_j^{\pm}(w)x_i^{\pm}(z), \\
& [x_i^{\pm}(z),x_j^{\pm}(w)]=0 \quad \text{ for $a_{ij}=0$}, \\
& x_i^{\pm}(z_1)x_i^{\pm}(z_2)x_j^{\pm}(w)-(q+q^{-1})x_i^{\pm}(z_1)
  x_j^{\pm}(w)x_i^{\pm}(z_2)+x_j^{\pm}(w)x_i^{\pm}(z_1)x_i^{\pm}(z_2) \\
& +\{ z_1\leftrightarrow z_2\}=0, \quad \text{for $a_{ij}=-1$}
\label{def}
\end{align*}
where
\begin{align} \delta(z)=\sum_{k\in \nz}z^k, \qquad {\rm and}\quad
   g_{ij}(z)=\frac{q^{a_{ij}}z-1}{z-q^{a_{ij}}}\quad |z|<|q^{a_{ij}}|.
\end{align}

 The  Hopf algebra structure
of $U_q^D(\hat{ \frak{g}})$ is given
by the following formulae.

\noindent{\bf Coproduct $\Delta$}
\begin{align*}
\text{}& \quad \Delta(q^c)=q^c\otimes q^c, \\
\text{}& \quad \Delta(x_i^+(z))=x_i^+(z)\otimes 1+
            \vphi_i(zq^{\frac{c_1}{2}})\otimes x_i^+(zq^{c_1}), \\
\text{}& \quad \Delta(x_i^-(z))=1\otimes x_i^-(z)+
            x_i^-(zq^{c_2})\otimes \psi_i(zq^{\frac{c_2}{2}}), \\
\text{}& \quad \Delta(\vphi_i(z))=
            \vphi_i(zq^{-\frac{c_2}{2}})\otimes\vphi_i(zq^{\frac{c_1}{2}}), \\
\text{}& \quad \Delta(\psi_i(z))=
            \psi_i(zq^{\frac{c_2}{2}})\otimes\psi_i(zq^{-\frac{c_1}{2}}).
\label{cop}
\end{align*}
\noindent{\bf Counit $\vep$}
\begin{align*}
\vep(q^c)=1 & \quad
\vep(\vphi_i(z))=\vep(\psi_i(z))=1, 
\quad \vep(x_i^{\pm}(z))=0.
\label{coun}
\end{align*}
\noindent{\bf Antipode $\quad a$}
\begin{align*}
\text{}& \quad a(q^c)=q^{-c}, \\
\text{}& \quad a(x_i^+(z))=-\vphi_i(zq^{-\frac{c}{2}})^{-1}
                               x_i^+(zq^{-c}), \\
\text{}& \quad a(x_i^-(z))=-x_i^-(zq^{-c})
                               \psi_i(zq^{-\frac{c}{2}})^{-1}, \\
\text{}& \quad a(\vphi_i(z))=\vphi_i(z)^{-1}, \quad
\text{} \quad a(\psi_i(z))=\psi_i(z)^{-1}.
\label{anti}
\end{align*}
The algebra $U_q^D(\hat{ \frak{g}})$ is a topological algebra. It means
 that $U_q^D(\hat{ \frak{g}})$ consists of certain series of the polynomials
 over the generators of $U_q^D(\hat{ \frak{g}})$ whose  action on the highest
 weight modules of $U_q(\hat{ \frak{g}})$ is well defined and the topology is
 a formal topology defined by natural filtration on the series.
 One of the possibilities to describe precisely such a topology is given
in \cite{KT2}. It is based on the construction of a Cartan-Weyl basis for
$U_q(\hat{ \frak{g}})$.

Let $e_{\pm\gamma}$ be the generators of Cartan-Weyl
 basis for a given normal ordering of positive roots (see details for affine
 case in \cite{KT1}), ${\frak{h}}$ be Cartan subalgebra of
 $U_q(\hat{ \frak{g}})$ generated by the elements
$k_{\alpha_{i}}^{\pm 1}$, $(i=0,1,\ldots,r)$. We consider
 formal  series on the following monomials
\begin{align}
e_{-\beta}^{n_{\beta}} \cdots e_{-\gamma}^{n_{\gamma}}
e_{-\alpha}^{n_{\alpha}}  e_{\alpha}^{m_{\alpha}}
e_{\gamma}^{m_{\gamma}} \cdots e_{\beta}^{m_{\beta}}
\label{CW16}
\end{align}
with coefficients from $U_{q}({\frak{h}})$, where
$\alpha<\gamma<\cdots<\beta$ in a sense of the fixed normal ordering
in the  system $\underline{\Delta}_{+}$ of positive roots with the condition
 that for any weight  $\lambda\in {\frak{h}}^*$ and for any
 positive integer $N$ there are only finitely many terms of total
 weight $\lambda$ satisfying the
condition
$$
\sum_{i=0}^{r}\sum_{\alpha \in \underline{\Delta}_{+}}(n_{\alpha}+
m_{\alpha})c_{i}^{(\alpha)} \leq N,
$$
where $c_{i}^{(\alpha)}$ are coefficients in a decomposition of the root
$\alpha$ with respect to the system  of simple roots of $\hat{ \frak{g}}$.
 The series  over the monomials
 (\ref{CW16}) with  additional constraint
\begin{align}
\sum_{i=0}^{r}\sum_{\alpha \in \Delta_{+}}
(n_{\alpha}+m_{\alpha})c_{i}^{(\alpha)}
> l  .
\label{CW18}
\end{align}
form a system of basic open sets in the algebra $U_q^D(\hat{ \frak{g}})$.

In an analogous manner,  the topology is introduced into the tensor square
 of $U_q^D(\hat{ \frak{g}})$.
Here we consider the series over monomials
\begin{align}
e_{-\beta}^{n_{\beta}} \cdots e_{-\gamma}^{n_{\gamma}}
e_{-\alpha}^{n_{\alpha}} e_{\alpha}^{m_{\alpha}}
e_{\gamma}^{m_{\gamma}} \cdots e_{\beta}^{m_{\beta}} \otimes
e_{-\beta}^{n'_{\beta}} \cdots e_{-\gamma}^{n'_{\gamma}}
e_{-\alpha}^{n'_{\alpha}} e_{\alpha}^{m'_{\alpha}}
e_{\gamma}^{m'_{\gamma}} \cdots e_{\beta}^{m'_{\beta}}
\label{CW19}
\end{align}
$\alpha<\gamma<\cdots<\beta$ with coefficients from
 $U_{q}({\frak{h}})\otimes U_{q}({\frak{h}}) $,
 such that for any two weights $\lambda, \mu \in
{\frak{h}}^*$ and for any positive interger $N$  , there are only finitely
 many terms of total weight $\lambda\otimes\mu$  satisfying the condition
$$
\sum_{i=0}^{r}\sum_{\alpha \in \Delta_{+}}(n_{\alpha}+n'_{\alpha}+
m_{\alpha}+m'_{\alpha})c_{i}^{(\alpha)}
\leq N
$$
The series  over the monomials
 (\ref{CW19}) with  additional constraint
\begin{align}
\sum_{i=0}^{r}\sum_{\alpha \in \Delta_{+}}
(n_{\alpha}+m_{\alpha} +n'_{\alpha}+ m'_{\alpha})c_{i}^{(\alpha)}
> l \ .
\label{CW22}
\end{align}
form a system of basic open sets in the algebra
 $U_q^D(\hat{ \frak g})\otimes U_q^D(\hat{ \frak g })$.

Another way to describe the topological structure of $U_q^D(\hat{ \frak{g}})$
 is to use instead of Cartan-Weyl basis of $U_q(\hat{ \frak{g}})$ the
 current type base introduced in \cite{DK} or the crystal base by Kashiwara.
 For \UqdvaD the Drinfeld generators are sufficient
 for the description of topology and of corresponding completion.

The  Hopf algebra structure   of $U_q^D(\hat{ \frak{g}})$ can be obtained
 from the Hopf structure of quantum affine algebra $U_q(\hat{ \frak{g}})$
 by means of twisting of the Hopf structure \cite{KT1}. Let us briefly
remind the  general idea of twistings \cite{R}.

Let ${\frak H}_{{\frak A}}:=({\frak A}, \Delta, S, \varepsilon, {\frak R})$
be a quasi-triangular
Hopf algebra with a comultiplication $\Delta$, an antipode $S$,
 counit $\epsilon$ and  universal R-matrix
${\frak R}$. Let also $F$, $F= \sum_i f_i\otimes f^i$ be an
invertible  element of
some extension $T ({\frak A }\otimes {\frak A})$ of ${\frak A}\otimes {\frak
A}$, such that the formula $$ \Delta^{(F)}(a):= F^{-1}\Delta(a)F,  \qquad a
\in {\frak A}, $$ determine a new comultiplication, i.e. $\Delta^{(F)}$
satisfies the coassociativity property
\begin{align} (\Delta^{(F)}\otimes {\rm
id})\Delta^{(F)}= ({\rm id}\otimes\Delta^{(F)})\Delta^{(F)}.  \label{b}
\end{align}
Then the comultiplication $\Delta^{(F)}$ is called the twisted
coproduct. This comultiplication defines a new quasitriangular
 Hopf algebra
${\frak H}_{{\frak A}}^{(F)}:=({\frak A}, \Delta^{(F)},
 S^{(F)}, \varepsilon,
{\frak R}^{(F)}).$
 Here $S^{(F)}(a)=u^{-1}S(a)u$, where
 $u:= (({\rm id} \otimes S)F) \circ 1 =
\sum_{i} f_{i} S(f^{i})$, and
$$
 {\frak  R}^{F}= (F^{21})^{-1}{\frak R}F.
$$
The condition (\ref{b}) is automatically satisfied if tensor $F$
 satisfies the  equalities:
\begin{align}
(\varepsilon\otimes {\rm id})F &=({\rm id}\otimes\varepsilon)F=1,\nonumber\\
({\rm id} \otimes F)(\Delta\otimes {\rm id})F &=
(F\otimes {\rm id})({\rm id} \otimes \Delta)F ,
\nonumber
\end{align}
One can also twist the  coalgebraic structure of the quasi-triangular
Hopf algebra ${\frak H}_{{\frak A}}$ =
$({\frak A}, \; \Delta, S, \epsilon,$ $ {\frak R})$
by using any automorphism in the algebraic sector ${\frak A}$.
Namely, let $\omega: {\frak A} \mapsto {\frak A}$ be an algebra
automorphism of ${\frak A}$. The  maps
 $\Delta^{(\omega)}: {\frak A} \mapsto {\frak A} \otimes {\frak A}$ and
$S^{(\omega)}: {\frak A} \mapsto {\frak A}$,
$$
\Delta^{(\omega)}(a) := (\omega \otimes \omega)\Delta(\omega^{-1}a) \quad
 S^{(\omega)}(a) := \omega S (\omega ^{-1}a) \quad a \in {\frak A}
$$
 define the quasitriangular Hopf algebra
 ${\frak H}_{\frak A}^{(\omega)}:=
({\frak A} \ , \Delta^{(\omega)}, S^{(\omega)}, \varepsilon,
{\frak R}^{(\omega)})$ where
$$
{\frak R}^{(\omega)} = (\omega \otimes \omega){\frak R} \ .
$$

Let ${\frak g}$ be any contragredient Lie algebra of finite
growth with a symmetrizable Cartan  matrix $A$. Denote by
$\underline{\Delta}_{+}$  the reduced system of  positive roots with
${\frak g}$, $\alpha_1, \ldots \alpha_r$ being positive
simple roots and  $\delta$ being the minimal positive imaginary
 root.
For the quantum algebra  $U_{q}({\frak g})$ with the coproduct rule as
$$\Delta(e_{\alpha_i})= e_{\alpha_i}\otimes 1+k_{\alpha_i}
\otimes e_{\alpha_i},$$
$$\Delta(e_{-\alpha_i})= 1\otimes e_{-\alpha_i}+
e_{-\alpha_i}\otimes k_{\alpha_i}^{-1},$$
we  have a  family of twisting, which can be described in both
 languages given above.

For a fixed normal ordering in $\underline{\Delta}_{+}$ and for any
$\gamma \, \in \underline{\Delta}_{+}$, $\gamma <\delta$ we put
\begin{align}
F_{\gamma} := \overrightarrow{\prod}_{\beta < \gamma} R^{21}_{\alpha}
\label{F}
\end{align}
where
$R_{\beta}= \exp_{q^{(\beta, \beta)}}(c(\beta)e_\beta\otimes e_{-\beta})$;
 $e_{\pm\beta}$ are Cartan-Weyl generators which correspond to this ordering
of the root system; $c(\beta)$ is a constant which is equal to $q^{-1}-q$ in
simple laced case in which we are interested in
 \cite{KT1}.
It is well known  \cite{KT1}, \cite{LS} that the
 tensor $F_{\gamma}$  defines a twisting of the Hopf structure
for $U_{q}({\frak g})$ (that is, the coassociativity condition
 (\ref{b}) is satisfied).

Let now $\gamma\in\underline{\Delta}_+$ is such that the initial segment
 $\gamma_1,\ldots ,\gamma_n = \gamma$ of the normal ordering is finite.
 Then $\gamma$ uniquely defines the element $\omega_\gamma$ of the Weyl
 group of ${\frak g}$. The element $\omega_\gamma$ together with its
 reduced decomposition $\omega_\gamma=s_{\alpha_1}\cdots s_{\alpha_n}$
can be defined by induction:
$$\omega_{\gamma_1}=id,\quad \omega_{\gamma_{k+1}}=\omega_{\gamma_k}
s_{\alpha_{i_{k}}}\Leftrightarrow \omega_{\gamma_k}(\alpha_{i_{k}})=
 \gamma_{k}.$$
 and uniquely defines the  automorphism $\widehat{\omega}_\gamma$ of
$U_{q}({\frak g})$ as
$$\widehat{\omega}_\gamma=T^{-1}_{\alpha_1}\cdots T^{-1}_{\alpha_n}$$
 where $T^{}_{\alpha_k}$ are Lusztig automorphisms. It is also known
 \cite{LS}, \cite{KT1} that
the twisting by $F_\gamma$ coincides with twisting by Lusztig automorphism
 $\widehat{\omega}_\gamma$ \cite{L}.

We are specially interested in the twist by a $q$-Weyl group element
 $\widehat{t}_{\delta}^N$ for  affine quantum algebras. Here
${t}_{\delta}$ is a translation by imaginary root $\delta$:
$$ t_{\delta}(\pm\alpha_k)=\pm\alpha_k \pm\delta.$$
\begin{prop}\label{D} \cite{KT1}.
(i) Quantum affine algebra $U_q(\widehat{\frak g})$ twisted by a tensor
 $F_\delta$, where $\delta$ is minimal imaginary root, is isomorphic to
 Drinfeld's realization $U_q^D(\widehat{\frak g})$ of quantum affine
 algebra

\noindent
(ii) The twisting by $F_\delta$ is equivalent to a limit of twisting
 by Lusztig automorphism $\widehat{t}_{\delta}^N$, $N\to\infty$ in the
 topology of
$U_q^D(\hat{ \frak{g}})\otimes U_q^D(\hat{ \frak{g}})$
 described above.
\end{prop}
 In particular, this proposition allows to describe the universal $R$-matrix
 for $U_q^D(\hat{ \frak{g}})$.
The universal R-matrix for $U_q(\widehat{\frak g})$ admits natural
 triangular decomposition \cite{KT1}, \cite{KST}:
\begin{align}
{\frak R}={\frak R}_+{\frak R}_0{\frak R}_-
\label{decom}
\end{align}
with
$${\frak R}_+=\prod_{\beta<\delta}^{\to}R_\beta,\quad
{\frak R}_0=R_\delta K,\quad {\frak R}_-=K^{-1}
\prod_{\beta>\delta}^{\to}R_\beta K ,$$
where $K=q^{-t_i\otimes t^i}$,
$t_i$ and $t^i$ are dual bases of extended Cartan subalgebra and $R_\delta$
 is a tensor, corresponding to imaginary root vectors, see precise expression
 in \cite{KT1}, \cite{KST}, where one should change in notations $q$ to
 $q^{-1}$ in order to adapt them to the definition of \UqgD.
 The Proposition \ref{D} implies that the universal R-matrix
 ${\frak R}^D$ for Drinfeld's realization has a form
\begin{align}
{\frak R}^D={\frak R}_0{\frak R}_-{\frak R}_+^{21} .
\label{univer}
\end{align}
Another presentation of the universal $R$-matrix for \UqgD  in a
form of integral over current operators is given in
\cite{DK}.
\section{The singular $R$-matrix}
Let us compute the projections of the universal R-matrices
 $(\widehat{t}_{2\delta}\otimes\widehat{t}_{2\delta})^N{\frak R}$ and
 of ${\frak R}^D$ onto tensor products of fundamental representations
 of \Uqn.

We start from $sl_2$ case. We fix $q$ in a region $|q|<1$.
With the definition of the quantum algebras
from \cite{KT1}, the projection of $\rho(z_1)\otimes \rho(z_2)$
of the universal $R$-matrix  of \Uqdva,  where
$\rho$ is two-dimensional  representation of $U_q(\widehat{sl}_2)$ and
$z_1$ and $z_2$ are formal variables,
gives the following for $z=z_1/z_2$:
$$
R(z)=r(z)\left(\begin{array}{cccc}
1&0&0&0\\
0&\frac{q^{}(1-z)}{1-q^{2}z}& \frac{(1-q^{2})}{1-q^{2}z}&0\\
0&\frac{(1-q^{2})z}{1-q^{2}z}&\frac{q^{}(1-z)}{1-q^{2}z}&0\\
0&0&0&1
\end{array}\right)
$$
where
$$r(z)=q^{-\frac{1}{2}}\frac{(zq^2;q^4)^2}{(z;q^4)(zq^4;q^4)}.$$
Note that $r(z)$ has a pole at $z=1$.
The triangular decomposition
of the universal R-matrix yields the triangular decomposition of the
matrix $R(z)$ \cite{KST}:
$$\begin{array}{lll}
R(z)&=r(z)&\left(\begin{array}{cccc}1&0&0&0\\
0&1&\frac{(q^{-1}-q^{})}{1-z}&0\\
0&0&1&0\\
0&0&0&1
\end{array}\right)
\left(\begin{array}{cccc}1&0&0&0\\
0&\frac{q^{}(1-q^{-2}z)}{1-z}&0&0\\
0&0&\frac{q^{}(1-z)}{1-q^{2}z}&0\\
0&0&0&1
\end{array}\right)
\left(\begin{array}{cccc}1&0&0&0\\
0&1&0&0\\
0&\frac{(q^{-1}-q^{})z}{1-z}&1&0\\
0&0&0&1
\end{array}\right)\end{array}$$ 
The tensor $F_N$ of the type (\ref{F}) defining the twist by
 $\widehat{t}^{2N}_{\delta}$, is equal to
$$\prod_{0\leq k\leq 2N-1}^{\to}\exp_{q^{2}}(q^{-1}-q^{})e_{\alpha+k\delta}
\otimes e_{-\alpha-k\delta}=
\prod_{0\leq k\leq 2N-1}^{\to}\exp_{q^{2}}(q^{-1}-q^{})x_k^+
\otimes x_k^-\ .$$
Its projection is a triangular matrix  $G(z)$
$$G(z)=\left(\begin{array}{cccc}
1&0&0&0\\ 0&1&\frac{(q^{-1}-q^{})(1-z^{2N})}{1-z}&0\\
0&0&1&0\\0&0&0&1\end{array}\right) .$$

Thus
 the image $R_{2N}(z)$ of twisted universal R-matrix
${\frak R}^{\tilde t^{2N}_{\delta}}$
  can be written as

\begin{align}
R_{2N}(z)=r(z)&\left(\begin{array}{cccc}1&0&0&0\\
0&1&\frac{(q^{-1}-q^{})z^{2N}}{1-z}&0\\
0&0&1&0\\
0&0&0&1
\end{array}\right)
\left(\begin{array}{cccc}1&0&0&0\\
0&\frac{q^{}(1-q^{-2}z)}{1-z}&0&0\\
0&0&\frac{q^{}(1-z)}{1-q^{2}z}&0\\
0&0&0&1
\end{array}\right)\times \nonumber\\
&
\left(\begin{array}{cccc}1&0&0&0\\
0&1&0&0\\
0&\frac{(q^{-1}-q^{})z^{-{2N}+1}}{(1-z)}&1&0\\
0&0&0&1
\end{array}\right)\nonumber\\
 =
q^{-1\over 2}\frac{(zq^2;q^4)^2}{(zq^4;q^4)^2}&\left(\begin{array}{cccc}
1/(1-z)&0&0&0\\
0&\frac{q^{}}{1-q^{2}z}& \frac{(1-q^{2})z^{2N}}{(1-q^{2}z)(1-z)}&0\\
0&\frac{(1-q^{2})z^{-{2N}+1}}{(1-q^{2}z)(1-z)}&\frac{q^{}}{(1-q^{2}z)}&0\\
0&0&0&1/(1-z)
\end{array}\right).
\end{align}
\begin{rem}
We would like to take certain the limit of the expression  above, for which we
need to define the topology on the space
End$(\rho\otimes \rho)\otimes [z , z^{-1}]$.
We define the topology by defing the limit of a series of
formal power series expression as the expression with each coefficient of
$z^n$ as   the limit of the each
 coefficient of   $z^n$ of the series. Clearly, this is the  only reasonable
topology we can have here. However, then,
there appears a subtle point about taking  the
the  limit of the operator $R_{2N}(z)$.
 There are two ways to take the limit
namely  take the limit of the expression above directly
 or take the limit of each component of the  Gauss decomposition
then multiply all the component together again.
These two limits are different.
The reason for such a disagreement comes from  the fact that, under the
topology we use, the following formula is not valid:
$$ \lim_{n\rightarrow \infty} f(n,z) \lim_{n\rightarrow \infty} g(n,z)
\neq \lim_{n\rightarrow \infty} f(n,z) g(n,z),$$
even both $ \lim_{n\rightarrow \infty} f(n,z)$
and $\lim_{n\rightarrow \infty} g(n,z)$ exist.
For example:
$$\lim_{n\rightarrow +\infty}z^n \Sigma_{m\in \Bbb Z_{+}} z^m= 0,\qquad
\lim_{n\rightarrow +\infty}z^{-n}
\Sigma_{m\in \Bbb Z_{+}} z^m= \delta(z),$$
but
$$\lim_{n\rightarrow +\infty}(z^n \Sigma_{m\in \Bbb Z_{+}} z^m)(
z^{-n} \Sigma_{m\in \Bbb Z_{+}} z^m)= (\Sigma_{m\in \Bbb Z_{+}} z^m)^2.$$
\end{rem}

Here,  we  would take the  second kind of limit namely the limit of the
product of the all the components of the Gauss decomposition, which
coincides with the topology we define in Section 1 for
affine quantum algebras.
Then we have
\begin{align}
{R^D(z)}=q^{-\frac{1}{ 2}}\frac{(zq^2;q^4)^2}{(zq^4;q^4)^2}
\left(\begin{array}{cccc}
1/(1-z)&0&0&0\\
0&\frac{q^{}(1-q^{-2}z)}{(1-z)^2}&0&0\\
0&\delta(z)&\frac{q^{}}{1-q^{2}z}&0\\
0&0&0&1/(1-z)
\end{array}\right).
\label{17}
\end{align}
We see, in accordance with Proposition \ref{D} that
$$R^D(z)=(R_+(z))^{-1}R(z)R_+(z)^{21}$$
 where $R_+(z)$ is the image of the first factor of
Gauss decomposition for $R(z)$;
 it coincides with the image of tensor $F_\delta$.
 So we prove
\begin{prop}
 The image $\rho(z_1)\otimes \rho(z_2) {\frak R}^{D}$ of the universal
 $R$-matrix for \UqdvaD is given by the formula (\ref{17}):
$${R^D(z)}=  \rho(z_1)\otimes \rho(z_2) {\frak R}^{D}. $$
\end{prop}

We also have
$$
((R^D)^{21})^{-1}(z)=r^{-1}(z)
\left(\begin{array}{cccc}
1-z&0&0&0\\
0&\frac{q^{-1}(1-z)^2}{1-q^{2}z}&(1-z)\delta(z)&0\\
0&0&{q^{}(1-q^{2}z)}&0\\
0&0&0&1-z
\end{array}\right)$$
which is diagonal.

The twisted operator ${\frak R}^D$ is well defined if we let it act on
$V_a\otimes V_b$, where $V_a$ and $V_b$ are the highest weight modules.
This operator may not well defined if we choose any
$V_a\otimes V_b$. The example above is  just a case that the
projection of ${\frak R}^D$ is well defined.

A similar  case is for \UqnD. Let
$V$ be the n-dimensional fundamental representation, the projection of
${\frak R}$, the universal R-matrix gives us:
\begin{align}
R(z) &=r(z)\left( \sum^n_{i=1}E_{ii}\otimes E_{ii} + \frac{z-1}{q^{}z-q^{-1}}
\sum^n\Sb i\neq
j\ \ i,j=1\endSb E_{ii}\otimes E_{jj}\right. \nonumber\\
&\left.- \frac{q-q^{-1}}{qz -q^{-1}}
\sum^n\Sb i>j\ \ i,j=1\endSb E_{ij}\otimes
E_{ij} -\frac{(q-q^{-1})z}{qz -q^{-1}} \sum^n\Sb i<j\ \
i,j=1\endSb E_{ij}\otimes E_{ji} \right),\nonumber
\end{align}
where $ r(z)=q^{-\frac{1}{ n}}
\frac{(zq^n;q^{2n})^2}{(z;q^{2n})(zq^{2n};q^{2n})}$
and $E_{ij}\in$ End$(V)$.

The triangular decomposition
of the universal $R$-matrix yields the triangular decomposition of the
matrix $R(z)$:
\begin{align}
R(z) &=r(z)\left( \sum\Sb i\neq j\endSb  E_{ii}\otimes E_{jj} +
\frac{(q-q^{-1})}{z-1}\sum^n\Sb i>j\ \ i,j=1\endSb E_{ij}\otimes
E_{ji}\right)\times \nonumber\\ &
\left(\sum^n\Sb i=1\endSb  E_{ii}\otimes E_{jj}+
\frac{q^{-1}z-q^{}}{z-1}\sum\Sb i<j\endSb  E_{ii}\otimes E_{jj}
+\frac{z-1}{q^{}z-q^{-1}}\sum\Sb i>j\endSb  E_{ii}\otimes E_{jj} \right)\times
\nonumber\\&
\left( \sum\Sb i\neq j\endSb  E_{ii}\otimes E_{jj}
 + \frac{(q-q^{-1})z}{z-1}\sum^n\Sb i<j\ \
i,j=1\endSb E_{ij}\otimes E_{ji} \right).
\nonumber
\end{align}

The tensor $F_N$  defining the twist by
 $\widehat{t}^{2N}_{\delta}$, is equal to
 ordered product of $R_\beta$ over real roots $\beta$ such that
 $t_{2N\delta}(\beta)\in -\underline{\Delta}_+$. These roots are
 $\varepsilon_j-\varepsilon_i+k\delta$, where  $i>j,\  0\leq k<2N(i-j)$.
Its projection is a triangular matrix $G_N(z)$
\begin{align}
G_N(z)=Id\otimes Id-\frac{(q-q^{-1})(1-z^{2N(i-j)})}{1-z}
\sum^n_{i>j,\, i,j=1} E_{ij}\otimes
E_{ji}
\label{gz}
\end{align}
 and
 the image $R_{2N}(z)$ of twisted universal R-matrix
${\frak R}^{\tilde t^N_{2\delta}}$ looks as

\begin{align}
R_{2N}(z)&= r(z)
\left( \sum\Sb \endSb  E_{ii}\otimes E_{jj} +z^{2N(i-j)}\frac{(q-q^{-1})}{z-1}
\sum^n\Sb i>j\ \ i,j=1\endSb E_{ij}\otimes
E_{ji}\right)\times\nonumber
\\ &\left(\sum^n\Sb i=1\endSb  E_{ii}\otimes E_{ii}+
\frac{q^{-1}z-q^{}}{z-1}\sum\Sb i<j\endSb  E_{ii}\otimes E_{jj}
+\frac{z-1}{q^{}z-q^{-1}}\sum\Sb i>j\endSb  E_{ii}\otimes E_{jj} \right)\times
\nonumber\\&
\left( \sum\Sb \endSb  E_{ii}\otimes E_{jj} +
 z^{2N(i-j)+1} \frac{(q-q^{-1})}{z-1}\sum^n\Sb i<j\ \
i,j=1\endSb E_{ij}\otimes E_{ji} \right).
\label{R2N}
\end{align}

and, finally,

$$
R^D(z)=q^{-\frac{1}{ n}}\frac{(zq^n;q^{2n})^2}{(zq^{2n};q^{2n})^2}
\left( \sum^n_{i=1}\frac 1 {(1-z)}E_{ii}\otimes E_{ii} +
\frac{q^{-1}z-q^{}}{(z-1)^2}\sum^n\Sb
 i<j\endSb E_{ii}\otimes E_{jj} +\right.$$
\begin{align}
+ \left.\frac{1}{(q^{}z-q^{-1})}\sum^n\Sb i>j=1\endSb E_{ii}\otimes E_{jj}
+ \delta(z)\sum^n\Sb i>j\endSb E_{ij}\otimes
E_{ji}\right)
\label{singular}
\end{align}
Again, we have that
$$
{R^D(z)}=  \rho(z_1)\otimes \rho(z_2) {\frak R}^{D}.
$$
We summarize the calculations  in the following proposition:
\begin{prop}
(i) The image $R^D(z)$ of the universal R-matrix in tensor product
 of fundamental representations of $U_q^D(\widehat{sl}_n)$ is given
 by the formula (\ref{singular}).

\noindent
(ii)
The matrix $R^D(z)$ is a product of the
 limit   $(N\to\infty)$ of all the components of the Gauss decomposition of
 the R-matrices $R_{2N}(z)$.

\noindent
(iii)
 The R-matrices $R_{2N}(z)$  satisfy the Yang-Baxter
 equation. They are connected to the standard trigonometric
solutions $R(z)$ of the Yang-Baxter equation by means of the relations
\begin{align}
 R_{2N}(z)=G_N^{-1}(z)R(z)G_N^{21}(z)
\label{twist}
\end{align}
or
\begin{align}
R_{2N}(z)=U^{2N}(z)R(z)U^{-2N}(z).
\label{conjugation}
\end{align}
where $G_N(z)$ is given by the relation (\ref{gz})
and $U(z)$ be the following diagonal matrix:
\begin{align}
U(z)=\sum_{i,j=1}^n z^{i-j}E_{ii}\otimes E_{jj}
\label{diagonal}
\end{align}
\end{prop}
One may treat the statement (iii) of the theorem as an existence of
 one-parameter family of symmetries for the trigonometric Yang-Baxter
 equation.

However this singular R-matrix, in some sense,  should not be called R-matrix

\begin{prop}
The matrix $R^D(z)$ does not satisfies Yang-Baxter equation.
\end{prop}
It can be checked directly that both
${{ R}^D} _{12}(z){{ R}^D} _{13}(zw){{ R}^D} _{23}(w)
$ and  ${{ R}^D} _{23}(w){{ R}^D} _{13}(zw){{ R}^D} _{12}(z)$
diverge.
Beyond this, we also know that  Drinfeld comultiplication
is not really well defined on finite dimensional representations
\cite{DI2} as well as the topology
 defined in section 1 does not
work in the case of finite dimensional representations.
Thus, we need to deal with convergence problem
of the twisted quaistriangular structure if we want to use $L$-operator's
 approach which is based on the properties of finite-dimensional
representations.

\section{Quasitrigular Hopf algebra structure, FRTS realization and
current operators for non-simple roots}
 Let us first remind the  technique of introducing the $L$-operators from
 the universal $R$-matrix \cite{FRT}, \cite{FR}.
  We can introduce the parameter dependence of the universal $R$-matrix
 for the algebra \UqgD and to reformulate the quasitriangular Hopf structure
 in the same manner as for quantum affine algebra \Uqg \cite{FR}.
This is done by
 means of automorphisms $D_z$, which acts on Chevalley generators as
$$D_z(e_{\alpha_i})= z^{\delta_{i,0}} e_{\alpha_i},
D_z(f_{\alpha_i})=z^{-1\delta_{i,0}}f_{\alpha_i}.  $$
The  action of $D_z$ is equivalent to conjugation by
$z^d$, where $d$ is the grading element:
 $$[d, e_{\pm\alpha_i}]={\pm\delta_{i,0}} e_{\pm\alpha_i}. $$
  We  define the map $\Delta_z(a)= (D_z \otimes
id)\Delta(a)$ and $\Delta'_z(a)= (D_z \otimes id)\Delta'(a)$, where
$a\in U_q(\hat {\frak g})$ and $\Delta'$ denotes the opposite
comultiplication. Then we put
 $${{\frak R}^D} (z)= (D_z\otimes {\text id})
q^{d\otimes c +c\otimes d}{{\frak R}^D}, $$
where ${\frak R}^D$ is as defined in (\ref{univer})  and $c$ is the central
 element of \UqgD. We can use the tensor ${{\frak R}^D} (z)$ for the
description of  quasitriangular structure of \UqgD on the category of highest
 weight modules. It means that on the tensor category of the highest weight
representations of  $U_q(\hat{\frak  g})$, the operator
 ${{\frak R}^D} (z)$
$\in$ $U_q(\hat{\frak  g})\hat  \otimes
 U_q(\hat {\frak  g})\otimes
{\Bbb C} [[z]]$  satisfies
$$ {{\frak R}^D} (z) \Delta_z (a)= (D^{-1}_{q^{c_2}} \otimes
D^{-1}_{q^{c_1}}) \Delta'_z (a) {{\frak R}^D} (z),$$
$$ (\Delta \otimes I){{\frak R}^D} (z)=
{{\frak R}^D} _{13}(zq^{c_2}){{\frak R}^D} _{23}(z),$$
$$(I\otimes\Delta){{\frak R}^D} (z)=
{{\frak R}^D} _{13}(zq^{-c_2}){{\frak R}^D} _{12}(z),
$$
$${{\frak R}^D} _{12}(z){{\frak R}^D} _{13}(zq^{c_2}w){{\frak R}^D} _{23}(w)=
{{\frak R}^D} _{23}(w){{\frak R}^D} _{13}(zq^{-c_2}w){{\frak R}^D} _{12}(z). $$
Here $c_1=c\otimes 1$, $c_2=1\otimes c$.

It follows  from the fact that the projection of
any of the operators ${\frak R}^+$, ${\frak R} ^0$ and ${\frak R}^-$ are
finite expressions when
they act on any vector $v_1\otimes v_2$ for $v_1$ and $v_2$ in any two
highest weight representations $V_1$ and $V_2$.

Let now $V$ be a finite dimensional   representation of
$U_q(\hat{\frak g})$. Let us  assume that the image  of
 ${{\frak R}^D}(\frac z w)$ on $V\otimes V$, which we is
denoted by  $ R^D(\frac z w)$ is well defined.
Let
\begin{align}
{ L}^{-}(z)= ({\text id} \otimes
\pi_{V})( {{\frak R}^D}(z^{-1})^{-1}),
\qquad
{ L}^{+}(z)= ({\text id} \otimes
\pi_{V}) ({{\frak R}^D} _{21}(z)).
\label{Lp}
\end{align}
Then both operators are well defined ans
\UqgD  as an algebra
 is  generated  by operator entries of
 ${ L}^{+}(z)$ and
${ L}^{-}(z)$ for the case ${\frak g}$ is ${\frak sl}_n$.
However, in general the equalities
\begin{align}
  R^D({z}/{w}) {L}^{\pm}_1(z)
 {L}^{\pm}_2(w) &= {L}^{\pm}_2(w)
{L}^{\pm}_1(z) R^D({z}/{w}),\\
 R^D({wq^{-c}}/{z}) {L}^+_1(z) {L}^-_2(w) &=
 {L}^-_2(w) {L}^+_1(z)
 R^D({wq^{c}}/{z}),
\label{frst}
\end{align}
 are not valid, in particular, both sides of (19) diverge.
Therefore, they can not be used to
describe the commutation relations
 between the generators of the algebra \UqgD. One can see this at the simplest
 example of \UqdvaD. Nevertheless, the analysis of the twisting procedure
 from the previous section allows to restore the full set of commutation
 relations between the entries of $L$-operators (\ref{Lp})
 and to prove that they are the defining relations of the algebra \UqgD.
 In the following we show this for the algebra \UqnD. From
 now on, in this section, we will deal  with \UqnD.

Let ${\frak R}$ be the universal R-matrix for $\Uqn$,
and  ${\frak R}^{\tilde t^{2N}_{\delta}}$ be the twisted universal R-matrix.
Let $V$ be a $n$ dimensional   representation of  fundamental representation
of $U_q(\hat{\frak g})$. The  image  $  R_{2N}(\frac z w)$ of
${\frak R}^{\tilde t^{2N}_{\delta}}$
on $V\otimes V$,  is given by the relation (\ref{R2N}).
Let $${ L}^-_{(2N)}(z)= ({\text id} \otimes
\pi_{V})( {{\frak R}^{\tilde t^{2N}_{\delta}}}(z^{-1})^{-1}), \qquad
{ L}^+_{(2N)}(z)= ({\text id} \otimes
\pi_{V}) ({{\frak R}^{\tilde t^{2N}_{\delta}}} _{21}(z)). $$
Both operators are well defined.
 As it follows from the definitions and from Proposition 3, the
 $L$-operators ${ L}^\pm_{(2N)}(z)$ satisfy the relations:
\begin{align}
 R_{2N}({z}/{w}) \Bigl({L}^{\pm}_{( 2N)}\Bigr)_1(z)
 \Bigl({L}^{\pm}_{( 2N)}\Bigr)_2(w) &= \Bigl({L}^{\pm}_{( 2N)}\Bigr)_2(w)
\Bigl({L}^{\pm}_{( 2N)}\Bigr)_1(z) R_{2N}({z}/{w}),\label{frstn1}\\
 R_{2N}({wq^{-c}}/{z}) \Bigl({L}^{+}_{( 2N)}\Bigr)_1(z)
\Bigl({L}^{-}_{( 2N)}\Bigr)_2(w) &=
 \Bigl({L}^{-}_{( 2N)}\Bigr)_2(w) \Bigl({L}^{+}_{( 2N)}\Bigr)_1(z)
 R_{2N}({wq^{c}}/{z}).
\label{frstn}
\end{align}
Moreover, on the category of highest weight representations,
$$
\lim_{N\rightarrow \infty}{L}^{\pm}_{( 2N)} (z) = {L}^{\pm}(z).
$$
We know also from the previous section, that
$$\lim_{N\rightarrow \infty}  R_{2N}({z}/{w})= R^{D}(z). $$
but still both sides of (\ref{frstn1}), (\ref{frstn}) diverge.
 Nevertheless, their matrix coefficients in highest weight module converge,
 so we have:
\begin{align}
<v,   R^{D}({z}/{w}) {L}^{\pm }_1(z)
 {L}^{\pm }_2(w) u> &=<v,  {L}^{\pm }_2(w)
{L}^{\pm }_1(z)  R^{D}({z}/{w}) u>,\cr
<v,  R^{D}({wq^{-c}}/{z}) {L}^{+}_1(z) {L}^{-}_2(w)u> &=
 <v, {L}^{-}_2(w) {L}^{+}_1(z)
  R^{D}({wq^{c}}/{z}) u>.
\end{align}
Here $u\in U$ is a vector of highest weight module $U$, $v\in U^*$ and both
sides of the equalities  are treated as analytic functions and
$\bar R^{D}({z}/{w})$ is a diagonal matrix, whose diagonal entries
are the same as that of $ R^{D}({z}/{w})$.
 As a corollary, we have the following statement.
\begin{prop}
In the category of the highest weight representations,
$$(z-w) R^D({z}/{w}) {L}^{\pm}_1(z)
 {L}^{\pm}_2(w) = (z-w)R^D({z}/{w}) {L}^{\pm}_2(w)
{L}^{\pm}_1(z),
$$
$$(z-q^{c}w)(z-q^{-c}w) R^D({wq^{-c}}/{z}) {L}^{+}_1(z)
 {L}^{-}_2(w) =(z-q^{c}w)(z-q^{-c}w) R^D({wq^{c}}/{z}) {L}^{-}_2(w)
{L}^{+}_1(z),
$$
Here $((z-w) R({z}/{w}))$, $(z^2-q^{2c}w^2) R^D({wq^{\pm c}}/{z})$
  are diagonal and
$z$, $w$ are both formal variables.
\end{prop}

{\bf Proof}. First we know that the equality above is
true on the level of matrix coefficients as analytic
functions. The limit of (21) and (22) diverge, it   tells us that there is no
poles on both side of the equations above, thus they are equal.

\begin{rem}
For the case of $\hat {\frak sl}_2$ and $\hat {\frak sl}_3$ ,
one can show by calculation that (20) is actually  valid. It seems this
should be true for any $n$, which we do not know how to prove in a simple way.
\end{rem}

Now we come to the Gauss coordinates of the $L$-operators.
Let
$$
L^+(z) =\pmatrix k^{-}_1(z)&&&0\\{} & \ddots & {}&{}\\ &
&\ddots&\\0 && &  k^{-}_n(z)\endpmatrix
 \pmatrix 1&&&0\\ e^{}_{12}(z) &\ddots & \\ &\ddots &\ddots\\
e_{1,n}(z)& &e^{}_{n-1,n}(z)&1\endpmatrix,
$$
$$
L^{-}(z) = \pmatrix 1 & f_{2,1}(z)&&f_{n,1}(z)\\ &\ddots
&\ddots\\ & &\ddots & f_{n,n-1}(z)\\0 & & &1\endpmatrix
\pmatrix k^{+}_1(z)&&&0\\ & \ddots & \\ &
&\ddots\\0 && &  k^{+}_n(z)\endpmatrix ,
$$
that is,
\begin{align}
L^+(z)=\Bigl(\sum_{i=1}^n k_i^- E_{ii}\Bigr)
\Bigl(1+\sum_{i<j}e_{i,j}E_{ji}\Bigr),
\qquad
L^-(z)=
\Bigl(1+\sum_{i>j}f_{i,j}E_{ji}\Bigr)
\Bigl(\sum_{i=1}^n k_i^+ E_{ii}\Bigr),
\label{lij}
\end{align}
and let also
$$
L^+_{(2N)}(z)=
\Bigl(1+\sum_{i>j}\tilde{f}_{i,j}^-E_{ji}\Bigr)
\Bigl(\sum_{i=1}^n \tilde{k}_i^- E_{ii}\Bigr)
\Bigl(1+\sum_{i<j}\tilde{e}_{i,j}^-E_{ji}\Bigr),
$$
$$
L^-_{(2N)}(z)=
\Bigl(1+\sum_{i>j}\tilde{f}_{i,j}^+E_{ji}\Bigr)
\Bigl(\sum_{i=1}^n \tilde{k}_i^+ E_{ii}\Bigr)
\Bigl(1+\sum_{i<j}\tilde{e}_{i,j}^+E_{ji}\Bigr),
$$
 be Gauss decomposition of finitely twisted $L$-operators.
Let us introduce the currents
$$x_i^{+}(z)=(q^{-1}-q)^{-1}e_{i,i+1}(zq^i),
\qquad
x_i^{-}(z)=(q^{-1}-q)^{-1}f_{i+1,i}(zq^i),
$$
$$\psi_i(z)=k^-_{i+1}(zq^iq^{-c/2} )^{-1}(k_i^-(zq^iq^{-c/2})), \qquad
\phi(z)=(k^+_{i+1}(zq^iq^{-c/2})^{-1}(k_i^+(zq^iq^{-c/2})),$$
and analogously for the finitely twisted $L$-operators:
$$\tilde{x}_i^{+}(z)=(q^{-1}-q)^{-1}\tilde{e}_{i,i+1}^-(zq^i),
\qquad
\tilde{x}_i^{-}(z)=(q^{-1}-q)^{-1}\tilde{f}_{i+1,i}^+(zq^i),
$$
$$\tilde{\psi}_i(z)=\tilde{k}^-_{i+1}(zq^iq^{-c/2})^{-1}
(\tilde{k}_i^-(zq^iq^{-c/2})), \qquad
\tilde{\phi}(z)=(\tilde{k}^+_{i+1}(zq^iq^{-c/2})^{-1}
(\tilde{k}_i^+(zq^iq^{-c/2})).$$
\begin{thm}

{\rm (i)}  The commutation relations between the currents
$\tilde{x}_i^{\pm}(z)$, $\tilde{\psi}_i(z)$ and $\tilde{\phi}(z)$ encoded
 in the relations (\ref{frstn1}), (\ref{frstn})
have well defined limit where they
 coincide with
the relations on ${x}_i^{\pm}(z)$, ${\psi}_i(z)$ and ${\phi}(z)$ in the
 definition of the algebra \UqnD.

{\rm (ii)} The comultiplication structure of \UqnD (see section 1) coincides
 with the following comultiplication rule of the $L$-operators:
$$ \Delta' ( L^+(z))=     L^+(zq^{-c_2})\otimes
  L^+(z),\qquad \Delta' ( L^-(z))=     L^-(z)\otimes
  L^-(zq^{c_1})
$$
\end{thm}

The proof  follows  from the results in \cite{KT1}. There are no convergency
 problems since the topology of \UqnD is compatible with Gauss decomposition
 of the universal $R$-matrix.
However the direct calculations can be treated as another direct  proof of
the results in \cite{KT1}.

\begin{rem}
 One may notice that in \cite{DF}, we shift
 the positive half and the negative
half currents  coming from FRTS realizations,
which is not necessary here. The reason is that the operator
$K$, which is one of the factor inside the universal R-matrix
includes operator similar to $D_1$, which automatically shifts the
half current just as in \cite{DF}. Let us also point out that
the comultiplication formula for
$L^{\pm}$ makes the  the  Drinfeld comultiplication transparent.
\end{rem}

Clearly, we know now everything about the diagonal entries and the
ones right above and below the diagonal entries of $L^{\pm}(z)$.
The question to write down the complete
commutation relation of all the current operators of
$L^{\pm}(z)$ arises.

Let $i=0, 1,...,n-1$ be the nods of the Dynkin diagram of
$\hat {\frak sl}_n$.
Let $\alpha_i$ $i=1, ..., n-1$ be the simple roots of ${\frak sl}_n$
corresponding to the nodes $i= 1...,n-1$.
Let $s_{\beta}$ denotes the
Weyl group element corresponding to the root $\beta$ of ${\frak sl}_n$.
Let us fix the following reduced decomposition of
the longest element $\omega_0$ of the Weyl group $W=S_n$ of Lie algebra
${\frak sl}_n$:
\begin{align}
\omega_0= \Bigl(s_{\alpha_1}s_{\alpha_2}...
s_{\alpha_{n-1} }\Bigr) \Bigl(s_{\alpha_1}s_{\alpha_2}...
s_{\alpha_{n-2} }\Bigr)...\Bigl(s_{\alpha_1}s_{\alpha_2}\bigr)
s_{\alpha_{1}}.\label{I} \end{align}
\begin{thm}
The entries of the second factor of $L^+(z)$ and
the first  factor of $L^-(z)$ are the same as the
current operators for all roots of ${\frak sl}_n$ corresponding to
the decomposition (\ref{I})
up to the shifts of  spectral parameters defined in \cite{DK}.
\end{thm}

{\bf Proof}. First, for any highest weight module $V_{\lambda}$ and
$v\in  V_{\lambda}$, we have
$ L^{\pm 2N}(z)v =  L^{\pm }(z)v, $
when $N$ is big enough.
With this, we can prove the theorem  for the case of
${\frak sl}_3$  by using the commutation relations for
$L^{\pm 2N}(z)$ as in \cite{DF},
then we factor out $k^{\pm}(z)$ operators and
 take the limit.

Below, we  give the  list of the commutation relations
for the case of $U_q(\hat {\frak sl}_3)$.
Let $K^+_{i,i+1}(z )= k^+_{i+1}(z q^{-c/2})^{-1}k^+_i(z q^{-c/2})$,
$\ K^-_{i,i+1}(z q^{-c/2})= k^-_{i+1}(z q^{-c/2})^{-1}k^-_i(z q^{-c/2})$.
 The currents  $e_{i,j}(z)$  and $f_{i,j}(z)$ are as defined in (\ref{lij}).
 First, the currents, corresponding to simple roots of ${\frak sl}_3$,
 satisfy the
relations given  in definition of \UqgD, if we shift them
accordingly as in Theorem 6.
We  give the rest in the following:
\begin{align*}
&  e_{1,2}(z_1)e_{2,3}(z_2)- \frac {zq-wq^{-1}}{z-w}
e_{2,3}(z_2)e_{1,2}(z_1)= \delta(z_1/z_2)e_{1,3}(z_1)  , \\
& f_{3,2}(z_1)f_{2,1}(z_2)- \frac {zq-wq^{-1}}{z-w}
f_{2,1}(z_2)f_{3,2}(z_1)= \delta (z_1/z_2)f_{3,1}(z_1),
\end{align*}
\begin{align*}
& [e_{1,3}(z),f_{3,1}(w)]=\frac{1}{q-q^{-1}}
  \left\{ \delta(\frac z wq^{-c})K^-_{1,2}(wq^{\frac{1}{2}c})
K^-_{2,3}(wq^{\frac{1}{2}c})-
          \delta(\frac z wq^{c})K^+_{1,2}(zq^{\frac{1}{2}c})
 K^+_{2,3}(zq^{\frac{1}{2}c})\right\}, \\
&  (zq-wq^{-1})e_{1,3}(z)e_{1,3}(w)=
(zq^{-1}-wq)e_{1,3}(w)e_{1,3}(z)  , \\
& (zq^{-1}-wq)f_{3,1}(z)f_{3,1}(w)=
(zq^{-1}-wq)f_{3,1}(w)f_{3,1}^+(z)  , \\
& [e_{1,2}(z),f_{3,1}(w)]=  \delta(\frac z wq^{c})K^+_{1,2}(zq^{\frac{1}{2}c})
f_{3,2}(zq^c)  , \\
& [e_{2,3}(z),f_{3,1}(w)]=  \delta(\frac z wq^{-c})f_{2,1}(w)
K^-_{2,3}(wq^{\frac{1}{2}c})  , \\
& [f_{2,1}(z),e_{1,3}(w)]=  \delta(\frac w zq^{-c})e_{2,3}(zq^{-c})
K^-_{1,2}(zq^{\frac{1}{2}c})  , \\
& [f_{3,2}(z),e_{1,3}(w)]= \delta(\frac w z q^{-c})
K^+_{1,2}(wq^{\frac{1}{2}c})e_{1,2}(w)  , \\
& \frac {(zq^{-1}-zq)}{z-w}f_{2,1}(z)f_{3,1}(w)=f_{3,1}(w)f_{2,1}(z) , \\
 &\frac  {z-w} {zq-wq^{-1}} f_{3,1}(w)f_{3,2}(z)= f_{3,2}(z)f_{3,1}(w)  , \\
& e_{1,2}(z)e_{1,3}(w)=\frac {(zq^{-1}-zq)}{z-w}
e_{1,3}(w)e_{1,2}(z) , \\
&  e_{1,3}(w)e_{2,3}(z)=\frac  {z-w} {zq-wq^{-1}}
 e_{2,3}(z)e_{1,3}(w).
\end{align*}
The last four  relations imply the cubic Serre relation.
Using induction, for general n, one can derive all the commutation
relations for the entries of $L^{\pm}(z)$ of  formulas for
$U_q(\hat {\frak sl}_n)$.

\section{Discussion}

This paper is more to use the known results to
explain the situation related to the quasitriangular structure
for the Drinfeld realization and ask the proper questions.
The quansitriangular structure for the
can not be recover completely Drinfeld realization on the
category of finite dimensional
representation.  The difficulty comes from the fact that we have not
been able to  properly control the convergence problem in
such a situation. However we want to emphasize that such a problem is
not  just for this  case  but rather general when we study the current type of
realizations of quantized algebras. The  similar situation happens to other
cases, such as Yangian, the new elliptic algebras \cite{KLP} \cite{EF1}.
We will discuss the situation looking  at the twisting of the R-matrix on
the finite dimensional representations.

For instance, for deriving current realization of the double of the Yangian for
 ${\frak sl}_2$,
we start from the application of tensor product of two-dimensional
 representation to the universal $R$-matrix and this gives the answer
 $$
  R(u)=r(u)\left(\begin{array}{cccc}
1&0&0&0\\
0&\frac{u}{u-i\hbar}& \frac{-i\hbar}{u-i\hbar}&0\\
0&\frac{-i\hbar}{u-i\hbar}&\frac{u}{u-i\hbar}&0\\
0&0&0&1
\end{array}\right)
$$
with
$$r(u)=\frac{\Gamma\Bigl(1-\frac{u}{2i\hbar}\Bigr)
\Gamma\Bigl(-\frac{u}{2i\hbar}\Bigr)}
{\Gamma^2\Bigl(\frac{1}{2}-\frac{u}{2i\hbar}\Bigr)}$$
The same procedure gives
$$
{R}^D(u)=\frac{\Gamma^2(1-\frac{u}{2i\hbar})}{\Gamma^2({1\over 2}
-\frac{u}{2i\hbar})}
\left(\begin{array}{cccc}
1&0&0&0\\
0&\frac{u+i\hbar}{u^2}&0&0\\
0&\delta(u)&\frac{1}{u-i\hbar}&0\\
0&0&0&1
\end{array}\right)
$$
and repeat all the arguments presented in previous sections. Analogously for
 the scaled elliptic algebra \cite{KLP}. Note that for a final result we
 do need the universal $R$-matrix, which is not known in this case. What is
 important is that the scalar factor before $R$-matrix, determined from
 crossing and unitary conditions, has a pole when $u=0$.
The recent work by \cite{EF1} can also be put in the same frame work,
where a Gauss decomposition of dynamic elliptic R-matrix is given.
\begin{align}
\nonumber
R(\zeta',\zeta)  =  &
A(\zeta, \zeta')
\left( 1 + \frac{\theta(\hbar)
\theta(\zeta' - \zeta + \lambda + \gamma)}
{\theta(\zeta' - \zeta) \theta(\lambda+\gamma)}
E_{2,1} \otimes E_{1,2} \right)
\\ & \nonumber
\left(
E_{1,1}\otimes E_{1,1} + E_{2,2} \otimes E_{2,2}
 +
{{\theta(\zeta'-\zeta)}\over{\theta(\zeta'-\zeta+\hbar)}}
E_{1,1}\otimes E_{2,2} +
{{\theta(\zeta'-\zeta-\hbar)}\over{\theta(\zeta'-\zeta)}}
E_{2,2} \otimes E_{1,1}
\right) 
\\ & \nonumber
\left( 1 - \frac{\theta(\hbar)
 \theta(\zeta - \zeta' + \lambda + \gamma)}
{\theta(\zeta - \zeta') \theta(\lambda + \gamma)}
 E_{1,2} \otimes E_{2,1} \right).
\end{align}
Here the function $A(\zeta, \zeta')$ has a pole at $\zeta=\zeta'$\cite{EF1}.
Let $A(\zeta, \zeta')=\frac {\bar A (\zeta, \zeta')}
{(\zeta-\zeta')}$, where ${\bar A (\zeta, \zeta')} $
is neither 0 nor a pole at  $\zeta-\zeta'=0$.
Similar procedure would give us the singular R-matrix:
\begin{align}
 \nonumber  &
\bar A(\zeta, \zeta')
\left(
{{1}\over {(\zeta-\zeta')}}E_{1,1}\otimes E_{1,1} +{{1}\over {(\zeta-\zeta')}}
 E_{2,2} \otimes E_{2,2}  +
{{\theta(\zeta'-\zeta)}\over{\theta(\zeta'-\zeta+\hbar)(\zeta-\zeta')}}
E_{1,1}\otimes E_{2,2}
\right.
\\ & \nonumber \left. +
{{\theta(\zeta'-\zeta-\hbar)}\over{\theta(\zeta'-\zeta)(\zeta-\zeta')}}
E_{2,2} \otimes E_{1,1}-\delta(\zeta-\zeta')E_{1,2} \otimes E_{2,1}
\right) \nonumber
\end{align}
We can see that this  singular R-matrix  does not depend on the dynamic
parameter anymore.
This   idea can be further applied to other cases,
especially the structure related to Elliptic R-matrices to derive Drinfeld
type of realizations using Gauss decomposition.
The importance of the
 new interpretation of Drinfeld realization
 from the point view of the structure is  that it gives us a better
way to understand the meaning related to the Gauss decomposition
of the FRTS realizations \cite{DF} \cite{KLP} \cite{EF1} \cite{KT3} and
the related quasi-triangular Hopf algebra structure for  Drinfeld type of
realizations.  Further more, we might use our idea to
to study the algebra
\cite{DI1}, namely to  understand those  algebra's
quasitriangular structure,
which is one of the motivations of this paper.

The entries of the $L^{\pm}(z)$ are also of great interest.
 From the calculation, we can see that the operator $f_1(z)$ and
$e_1(z)$ commute with themselves for the case of
$U_q(\hat{\frak sl}_n)$. This operator actually is used
in the quantum semi-infinite construction \cite{DFe1}.
For the case of
$U_q(\hat{\frak sl}_n)$, the matrix operator $L^{\pm} (z)$ are triangular,
whose entries should give us naturally the current operators that
correspond to other non-simple roots. Those  current operator will
enable us to derive a better description of the integrable conditions
\cite{DM} as in the classical case and give the semi-infinite
construction $U_q(\hat{\frak sl}_n)$.
One more  possible application is that
the current operators corresponding to the longest root may naturally
provide us a way to derive representations of quantum toroidal algebras.

Another immediate application  is to use the idea  to
derive the intertwiners corresponding to the Drinfeld realizations
\cite{DI2} and the other way around. Namely, we can compose the
intertwiners with the twister $D$ to derive the intertwiners corresponding
to the Drinfeld comultiplication and vice versa. Combining
this with the results \cite{FR}, we can obtain
a   q-KZ equation corresponding to the Drinfeld comultiplication
 \cite{DFe2}.
This new interpretation will definitely help us to understand
problem related to the intertwiners
 corresponding to the Drinfeld comultiplication, for example,
why certain types of  the intertwiners
 corresponding to the Drinfeld comultiplication do not exist\cite{DI1}.

For the case of $U_q(\hat{\frak sl}_2)$  \cite{Ke},
a  FRTS realization using a diagonal R-matrix is given by a diagonal R-matrix.
There is no details about  the main results  Proposition 3.1
about the FRTS realization, which mathematically looks impossible.

\bigskip

{\it Acknowlegements.}\
The authors are thankful to Prof. B. Feigin for useful
discussions. S.Kh would like to thank Prof. T. Miwa and RIMS for hospitality
 in Kyoto, where the work was finished.
 S.Kh. was supported in part by RFBR grant
  98-01-00344, grant  96-15-96455  for support of scientific schools and
INTAS grant 930166-ext 
and by Award  No. RM2-150 of the U.S. Civilian Research \& Development
 Foundation (CRDF) for the Independent States of the Former Soviet Union.
J. D. would  like to  thank the
hospitality of the Feigins during his visit in Moscow, where this work
was started.

\end{document}